\documentclass[10pt]{IEEEtran}
\usepackage{epsfig,color,amsmath}
\usepackage[noadjust]{cite}
\usepackage{amsthm}
\IEEEoverridecommandlockouts
\usepackage{wrapfig} 
\usepackage{bm}
\usepackage{epstopdf}
\usepackage{amssymb}
\usepackage{url}
\usepackage{enumitem}
\usepackage{multirow}
\usepackage{makecell}
\usepackage{algorithm,algorithmic}	
\setlist[itemize]{leftmargin=*}

\makeatother

\newtheorem{theorem}{Theorem}
\newtheorem{proposition}{Proposition}

\newtheorem{mydef}{Definition}
\newtheorem{mylem}{Lemma}

\newtheorem{rem}{Remark}

\newtheorem{asmp}{Assumption}

\newcommand{\normof}[1]{\|#1\|}

\newcommand{\m}{\boldsymbol}
\allowdisplaybreaks[4]
\newcommand{\inp}[2]{\langle {#1},{#2} \rangle} 
\newcommand{\nrm}[1]{\vert\vert #1 \vert\vert}
\newcommand{\pder}[2]{\frac{\partial #1}{\partial #2}}

\usepackage{graphicx}
\usepackage{ifpdf}
\usepackage{epstopdf}
\usepackage{ifpdf}
\ifpdf
\DeclareGraphicsExtensions{.eps}
\else
\DeclareGraphicsExtensions{.eps}
\fi
\usepackage[caption=off]{subfig} 

\usepackage[colorlinks = true,
linkcolor = blue,
urlcolor  = blue,
citecolor = blue,
anchorcolor = blue]{hyperref}
\usepackage{stackengine} 
\stackMath
\usepackage[flushleft]{threeparttable}
\usepackage{makecell}

\usepackage{multicol}
\usepackage{mathtools}

\usepackage[normalem]{ulem}
\usepackage{mathrsfs}
\DeclareMathAlphabet\mathbfcal{OMS}{cmsy}{b}{n}
\DeclarePairedDelimiter\norm{\lVert}{\rVert}%

\usepackage{stackengine}

\newcommand{\ubar}[1]{\text{\b{$#1$}}}
\title{{Computing Lipschitz Constants for Hydraulic Models of Water Distribution Networks}}
\author{Ramsey Shadfa$\text{n}^\dagger$, Shen Wan$\text{g}^\dagger$, Sebastian A. Nugroh$\text{o}^\dagger$, Fengxin Che$\text{n}^\ddagger$, and Ahmad F. Tah$\text{a}^{\dagger,\ast}$
	\vspace{-0.64cm}
	\thanks{$^\dagger$Department of Electrical and Computer Engineering, The University of Texas at San Antonio, TX 78249. Emails: \{ramsey.shadfan, shen.wang,sebastian.nugroho,ahmad.taha\}@utsa.edu.}
	\thanks{$^\ddagger$Department of Mathematics, The University of Texas at San Antonio, TX 78249. Email: fengxin.chen@utsa.edu}
	\thanks{$^\ast$Corresponding author. }
	\thanks{This work is partially supported by Valero Energy Corporation and National Science Foundation under Grants 1728629 and 1917164.}}
\begin{document}

\maketitle

\setlength{\abovedisplayskip}{3.5pt}
\setlength{\belowdisplayskip}{3.5pt}
\setlength{\abovedisplayshortskip}{3.1pt}
\setlength{\belowdisplayshortskip}{3.1pt}

\newdimen\origiwspc%
\newdimen\origiwstr%
\origiwspc=\fontdimen2\font
\origiwstr=\fontdimen3\font

\fontdimen2\font=0.64ex

\begin{abstract}
Drinking water distribution networks (WDN) are large-scale, dynamic systems spanning large geographic areas. Water networks include various components such as junctions, reservoirs, tanks, pipes, pumps, and valves. Hydraulic models for these components  depicting mass and energy balance form nonlinear algebraic differential equations (NDAE). While control theoretic studies have been thoroughly explored for other complex infrastructure such as power and transportation systems, little is understood or even investigated for feedback control and state estimation problems for the NDAE models of WDN. The objective of this paper is to showcase a complete NDAE model of WDN followed by computing Lipschitz constants of the vector-valued nonlinearity in that model. The computation of Lipschitz constants of hydraulic models is crucial as it paves the way to apply a plethora of control-theoretic studies for water system applications. In particular, the computation of Lipschitz constant is explored through closed-form, analytical expressions as well as via numerical methods. Case studies reveal how such computations fare against each other for various water networks. 
\end{abstract}

\vspace{-0.34cm}
\begin{IEEEkeywords}
	Hydraulic Models, Water Distribution Networks, Nonlinear dynamic systems, Lipschitz Continuity,  One-Sided Lipschitz, Nonlinearity Parameterization. 
\end{IEEEkeywords}

\vspace{-0.34cm}

\section{Introduction and Paper Contributions}

\IEEEPARstart{I}{n} the past few decades, virtually thousands of control-theoretic methods have been proposed to perform  either real-time feedback control or dynamic state estimation in dynamic networked systems. While a significant number of these methods and studies are tailored for linear---or even linearized models around operating points---dynamic system models, many other studies investigated methods for the more general and realistic nonlinear models. The three-decade old textbook~\cite{boyd1994linear} provides a  listing of such literature with a focus on semidefinite programming-based formulations for mostly linear dynamic models.

As for feedback control and state estimation for nonlinear systems, the majority of studies assume that the involved nonlinear functions in the dynamic models are: \textit{Lipschitz continuous}~\cite{Yadegar2018,Phanomchoeng2010,ALESSANDRI2004,ichalal2012observer,zemouche2013lmi}, or even a less restrictive assumption such as \textit{one-sided Lipschtiz continuity}~\cite{Abbaszadeh2010,zhang2012full,benallouch2012observer}, or more restrictive such as the \textit{bounded Jacobian} assumption~\cite{Phanomchoeng2010b,Jin2018,Wang2018}; the actual definitions of some of these assumptions are given in the ensuing sections and in the references therein. A key ingredient to such studies is the computation of Lipschitz constants of the involved nonlinear models. Specifically, the design of controller or estimator is reliant on computing Lipschitz constants of the encompassed vector-valued nonlinearity.
Although a wide range of nonlinear, control-theoretic algorithms have been investigated for electric power network and transportation systems with focus on feedback control or state estimation applications, not much is understood or investigated in terms of understanding the nonlinear, dynamic functions modeling hydraulics in drinking water distribution networks (WDN).  

\noindent \textbf{{Paper Objective and Contributions}}  \hspace{0.1cm}  The objective of this paper is to investigate analytical and computational methods to compute and bound Lipschitz and one-sided Lipschitz constants for a generalized, realistic model of the hydraulics in WDN comprised of reservoirs, pumps, valves, and tanks.  In this paper, we compare the derived analytical results with the numerical ones for different WDN. We show that the vector-valued nonlinearity in the WDN dynamics has a special structure that allows closed-form analytical expressions of the Lipschitz constants. Applying such Lipschitz constant computation to feedback control and state estimation problem formulations~\cite{Yadegar2018,Phanomchoeng2010,ALESSANDRI2004,ichalal2012observer,zemouche2013lmi,Abbaszadeh2010,zhang2012full,benallouch2012observer} is beyond the scope of this paper, but is an important future research direction of the authors. 

In particular, this paper considers arbitrary WDN topologies, different types of valves, and assumptions that are satisfied in realistic WDN hydraulic simulations. The computation of such constants provides three important contributions to the literature: \textit{(i)} bounds on how the nonlinearities encompassed in hydraulic models grows and evolves; \textit{(ii)} the utilization of such constants to design state feedback controllers and state estimators; \textit{(iii)} a thorough investigation of analytical and computational values of Lipschitz constants for various water networks. The paper organization is given as follows. 

\noindent \textbf{{Paper Organization}}  \hspace{0.1cm} Section~\ref{sec:WDNModel} presents the general differential algebraic equation model of WDN hydraulics. The model essentially presents mass and energy balance in all {components}  and is given in both continuous- and discrete-time. Section \ref{sec:pbmform} presents the problem formulation while Sections~\ref{sec:Lipschitz} and~\ref{sec:OSL} present the theoretical contribution of the paper through computing Lipschitz and one-sided Lipschitz constants. Section~\ref{sec:global} presents a brief discussion on computing such constants using global optimization and numerical algorithms. Section~\ref{sec:tests} concludes the paper with numerical tests on various water distribution network topologies and parameters. The paper's broader impact is to encourage advanced control-theoretic algorithms for nonlinear models in WDN.

\noindent \textbf{{Paper Notation}}  \hspace{0.1cm} The symbols $\mathbb{R}^n$ and $\mathbb{R}^{p\times q}$ denote column vectors with $n$ elements and real-valued matrices with size $p$-by-$q$.  Italicized, boldface upper and lower case characters represent matrices and column vectors---$a$ is a scalar, $\m a$ is a vector, and $\m A$ is a matrix. Matrix $\m I_n$ is a $n\times n$ identity square matrix, while $\m 0$ and $\m O$ represent zero vectors and matrices of appropriate dimensions. 

\vspace{-0.32cm}
\section{Water Distribution Networks (WDN) Dynamics}~\label{sec:WDNModel}
\vspace{-0.32cm}

We model WDN by a directed graph $(\mathcal{W},\mathcal{E})$.  Set $\mathcal{W}$ defines the nodes and is partitioned as $\mathcal{W} = \mathcal{J} \bigcup \mathcal{T} \bigcup \mathcal{R}$ where $\mathcal{J}$, $\mathcal{T}$, and $\mathcal{R}$ stand for the collection of junctions, tanks, and reservoirs. Let $\mathcal{E} \subseteq \mathcal{W} \times \mathcal{W}$ be the set of links, and define the partition $\mathcal{E} = \mathcal{P} \bigcup \mathcal{M} \bigcup \mathcal{V}$, where $\mathcal{P}$, $\mathcal{M}$, and $\mathcal{V}$ stand for the collection of pipes, pumps, and valves. For the $i^\mathrm{th}$ node, set $\mathcal{N}_i$ collects its neighboring nodes and is partitioned as $\mathcal{N}_i = \mathcal{N}_i^\mathrm{in} \bigcup \mathcal{N}_i^\mathrm{out}$, where $\mathcal{N}_i^\mathrm{in}$ and $\mathcal{N}_i^\mathrm{out}$ stand for the collection of inflow and outflow nodes.  Tab.~\ref{table:sets} summarizes the set and variables notation used in this paper.
\begin{table}[t]
	\fontsize{7}{7}\selectfont
	\caption{Set and Variable notation.}
		\vspace{-0.2cm}	
	\centering
	\renewcommand{\arraystretch}{1.5}
	\begin{tabular}{ c|c }
			\hline 
		\textit{Notation} & \textit{Set/Variable Notation Description} \\ \hline
		$\mathcal{W}$& A set of nodes including junctions, tanks and reservoirs   \\ \hline
		\makecell{$\mathcal{E}$}
		& \makecell{A set of links including pipes, pumps and valves }  \\	\hline
		$\mathcal{J}, \mathcal{T}, \mathcal{R}$ &  A set of $n_j$ junctions, $n_t$ tanks, and $n_r$ reservoirs    \\ \hline
		$\mathcal{P}$, $\mathcal{M}$ , $\mathcal{V}$ &  A set of $n_p$ pipes, $n_m$ pumps, and $n_v$ valves   \\ \hline
		$\mathcal{N}_i$ &  A set of neighbors node of the $i^\mathrm{th}$ node, $i \in \mathcal{W}$  \\  \hline
		$d_i$, $h_i$&  Demand and head at node $i$   \\ \hline
		$h_i^{\mathrm{TK}}$, $h_i^{\mathrm{R}}$  &  Head at  tank $i$ and  reservoir $i$    \\ \hline
		$h_{ij}^{\mathrm{P}}$, $h_{ij}^{\mathrm{M}}$, $h_{ij}^{\mathrm{V}}$  &  Head loss for the pipe, pump, and valve from $i$ to $j$   \\ \hline
		\makecell{$q_{ij}$} & \makecell{Flow through a pipe, valve or pump from node $i$ to node $j$}  \\	\hline
		\hline			
	\end{tabular}
\vspace{-0.3cm}
	\label{table:sets}
\end{table}
According to these basic laws, the equations that model mass and energy conservation for all components in WDN can be written in explicit and compact matrix-vector forms in Tab.~\ref{tab:models}. 

\subsection{Models of Components}~\label{sec:Model_iass}

\vspace{-0.4cm}

\subsubsection{Tanks and Reservoirs}	The water volume dynamics in the $i^\mathrm{th}$ tank at time $k$ can be expressed by a discrete-time difference equation~\eqref{equ:tank-volume}, while the head created by the tank can be described as~\eqref{equ:tank-head}
\vspace{-0.4em}

\begin{subequations}
\small	\begin{align}
 \hspace{-12pt}	V_{i}^\mathrm{TK}(k+1) &= V_{i}^\mathrm{TK}(k) + \Delta t \hspace{-3pt}\left(\hspace{-1pt}\sum_{j \in \mathcal{N}_i^\mathrm{in}} \hspace{-2pt} q_{ji}(k)\hspace{-2pt}-\hspace{-6pt}\sum_{j \in \mathcal{N}_i^\mathrm{out}}\hspace{-4pt} q_{ij}(k)\hspace{-3pt}\right)\hspace{-3pt}  ~\label{equ:tank-volume} \\
	h_i^{\mathrm{TK}}(k) &= \frac{V_i(k)}{A_i^{\mathrm{TK}}} + \mathscr E_i^{\mathrm{TK}}, \; i\in\mathcal{T}, ~\label{equ:tank-head} 
	\end{align}
\end{subequations}
where $V_{i}^\mathrm{TK}$ and $\Delta t$ are the volume and sampling time; $q_{ji}(k),\;i \in \mathcal{J},\;j \in \mathcal{N}_i^\mathrm{in} $ stands for the inflow of the $j^\mathrm{th}$ neighbor, while $q_{ij}(k),\;i \in \mathcal{J},\;j \in \mathcal{N}_i^\mathrm{out} $ stands for the outflow of the $j^\mathrm{th}$ neighbor; $h_i^{\mathrm{TK}}$, $A_i^{\mathrm{TK}}$, and $\mathscr E_i^{\mathrm{TK}}$ respectively stand for the head, cross-sectional area, and elevation of the $i^\mathrm{th}$ tank.
Combining~\eqref{equ:tank-volume} and~\eqref{equ:tank-head}, the head changes from time $k$ to $k+1$ of the $i^\mathrm{th}$ tank can be written as~\eqref{equ:tankhead} in Tab.~\ref{tab:models}. 
We assume that reservoirs have infinite water supply and the head of the $i^\mathrm{th}$ reservoir is fixed~\cite[Chapter 3.1]{rossman2000epanet},~\cite[Chapter 3]{puig2017real} and presented as~\eqref{equ:head-reservoir} in Tab.~\ref{tab:models}.


\subsubsection{Junctions and Pipes}
Junctions are the points where water flow merges or splits. The expression of mass conservation of the $i^\mathrm{th}$ junction at time $k$ can be written as~\eqref{equ:nodes} in Tab.~\ref{tab:models}, and $d_i$ stands for the demand that is extracted from node $i$.
%

The major head loss of a pipe from node $i$ to $j$ is due to friction and is determined by~\eqref{equ:head-flow-pipe} from Tab.~\ref{tab:models}, where $R$ is the resistance coefficient and $\mu$ is the constant flow exponent in the Hazen-Williams formula. Minor head losses are ignored in this paper.


\subsubsection{Head Gain in Pumps} 	A head increase/gain can be generated by a pump between the suction node $i$ and the delivery node $j$. The pump properties decide the relationship function between the pump flow and head increase~\cite{linsley1979water}, \cite[Chapter 3]{rossman2000epanet}. Generally, the head gain can be expressed as~\eqref{equ:head-flow-pump}, where $h^\mathrm{s}$ is the shutoff head for the pump; $q_{ij}$ is the flow through a pump; $s_{ij} \in [0,s_{ij}^{\mathrm{max}} ]$ is the relative speed of the same pump; $r$ and $\nu$ are the pump curve coefficients.

\subsubsection{Valves} In this work, we focus on general purpose valves (GPV). 
{GPVs} are used to represent a link with a special flow-head loss relationship instead of one of the standard hydraulic formulas. They can be used to model turbines, well draw-down or reduced-flow backflow prevention valves~\cite[Chapter 3.1]{rossman2000epanet}. In this paper, we assume that the GPVs are modeled as a pipe with controlled resistance coefficient and can be expressed as~\eqref{equ:head-flow-valve} in Tab.~\ref{tab:models}, where $o_{ij} \in [0,1]$  is a control variable depicting the openness of a valve. The model can be seamlessly extended to other types of valves such as pressure reducing of flow control valves.

\begin{table}[t]
	\fontsize{7}{7}\selectfont
	\centering
	\renewcommand{\arraystretch}{1.2}
	\caption{{Detailed hydraulic model in WDN}. }
	\vspace{-0.2cm}			
	\begin{tabular}{c|c}
		\hline 
		{\textit{Types}} 	& {\textit{Original Hydraulic Model}} \\ \hline
		{\textit{Tanks}} 
		&
		
		\parbox{7cm}{
			\begin{align}~\label{equ:tankhead}
			\hspace{-10pt}h_{i}^{\mathrm{TK}}(k\hspace{-2pt}+\hspace{-2pt}1) \hspace{-2pt}=\hspace{-2pt} h_{i}^{\mathrm{TK}}(k) \hspace{-2pt}+\hspace{-2pt} \frac{\Delta t}{A_i^{\mathrm{TK}}}\left(\hspace{-1pt}\sum_{j \in \mathcal{N}_i^\mathrm{in}}\hspace{-3pt}q_{ji}(k)\hspace{-2pt}-\hspace{-6pt}\sum_{j \in \mathcal{N}_i^\mathrm{out}} \hspace{-3pt}q_{ij}(k)\hspace{-3pt}\right)
			\end{align}	\vspace{-0.2em}
		}
		\\
		\hline
		{\textit{Reservoirs}} 
		& 	\parbox{7cm}{	\vspace{-0.2em}
			\begin{align}~\label{equ:head-reservoir}
			h_i^{\mathrm{R}}(k) = h_i^{\mathrm{R}}
			\end{align}
			\vspace{-0.9em}
		}
		\\
		\hline
		{\textit{\makecell{Junction\\ nodes}}} 
		&
		\parbox{7cm}{
			\vspace{-0.2em}
			\begin{align}~\label{equ:nodes}
			\sum_{j \in \mathcal{N}_i^\mathrm{in}} q_{ji}(k) - \sum_{j \in \mathcal{N}_i^\mathrm{out}} q_{ij}(k) = d_i(k)
			\end{align}
			\vspace{-0.5em}
		}
		\\
		\hline

		{\textit{Pipes}}
		&  \parbox{7cm}{
			\vspace{-0.3em}
			\begin{align}~\label{equ:head-flow-pipe}
			h_{ij}^\mathrm{P}(k)  = h_{i}(k) - h_{j}(k) = R_{ij} {q_{ij}(k)}|q_{ij}(k)|^{\mu-1}
			\end{align}
			\vspace{-0.9em}
		}
		\\
		\hline

		{\textit{Pumps}} 
		& \parbox{7cm}{
			\vspace{-0.3em}
			\begin{align} \label{equ:head-flow-pump}
			\hspace{-7pt} h_{ij}^\mathrm{\mathrm{M}}(k) = h_{i}(k) - h_{j}(k) = -{s_{ij}^2(k)}(h^\mathrm{s} - r_{ij}  (q_{ij} s_{ij}^{-1})^\nu )
			\end{align}
			\vspace{-0.8em}
		} 
		\\
		\hline
		{\textit{\hspace{-5pt}Valves \hspace{-5pt}}} 
		&
		
		\parbox{7cm}{
			\vspace{-0.3em}
			\begin{align}~\label{equ:head-flow-valve}
			\hspace{-9pt}h_{ij}^\mathrm{V}(k)  = h_{i}(k) - h_{j}(k) = o_{ij}(k) R_{ij} {q_{ij}(k)}|q_{ij}(k)|^{\mu-1}
			\end{align}
			\vspace{-0.95em}
		}
		\\
		%
		%
		%
		\hline \hline
	\end{tabular}
	\label{tab:models}%
	\vspace{-1.5em}
\end{table}
\vspace{-0.3cm}


	
\subsection{Control-oriented WDN model}
The WDN model in the previous section can be abstracted to difference algebraic equation (DAE) as~\eqref{equ:dae-abstract}. Define $\m x_1$, $\m x_2$, and  $\m x_3$ as vectors collecting heads at $n_j$ junctions, $n_r$ reservoirs, and $n_t$ tanks; define $\m u$ and $\m v$ as the vectors collecting flows through controllable and uncontrollable elements, e.g., $n_v+n_m$  pumps and valves belong to controllable set, and $n_p$ pipes are uncontrollable. Note that we assume all valves are GPVs. Collecting the mass and energy balance equations of tanks~\eqref{equ:tankhead}, reservoirs~\eqref{equ:head-reservoir}, junctions \eqref{equ:nodes}, pipes~\eqref{equ:head-flow-pipe}, pumps~\eqref{equ:head-flow-pump}, and valves \eqref{equ:head-flow-valve}, we obtain the following control-oriented DAE model
\begin{subequations}~\label{equ:dae-abstract}
	\begin{align}
	\hspace{-3em}\m x_3(k + 1) &\hspace{-1pt}=\hspace{-1pt}  \m x_3(k) + \m B_v \m v(k) ~\label{equ:tankhead-abstracted} \\
	\m 0_{n_p} &\hspace{-1pt}=\hspace{-1pt}  \m E_{x_1} \m x_1(k) \hspace{-2pt}+\hspace{-2pt} \m E_{x_2} \m x_2(k) \hspace{-2pt}+\hspace{-2pt} \m E_{x_3} \m x_3(k) \hspace{-2pt}+\hspace{-2pt}\m f^\mathrm{P}(\m v)  ~\label{equ:pipeloss-abstracted}\\
	\m 0_{n_m+n_v} &\hspace{-1pt}=\hspace{-1pt} \m F_{x_1} \m x_1(k) + \m F_{x_2} \m x_2(k)  
	{}+ \m f^\mathrm{MV}(\m u)    ~\label{equ:pumpgain-abstracted}\\
	\m 0_{n_j} &\hspace{-1pt}=\hspace{-1pt} \m G_{u} \m u(k) + \m G_{v} \m v(k) +\m d(k)  ~\label{equ:demand-abstracted} \\
	\m 0_{n_r} &\hspace{-1pt}=\hspace{-1pt} \m H_{x_2} \m x_2(k)  +  \m h^{\mathrm{R}} ~\label{equ:reservoir-abstracted} 
	\end{align}
\end{subequations}
where $\m B_{\bullet}$, $\m E_{\bullet}$, $\m F_{\bullet}$, $\m G_{\bullet}$, and  $\m H_{\bullet}$ are constant matrices that depend on the WDN topology and the aforementioned hydraulics and $\m 0_{n}$ is a zero-vector of size $n$. 
 Let $\m f^\mathrm{P}(\m v), \m f^\mathrm{M}(\m u),$ and $\m f^\mathrm{V} (\m u)$ be the nonlinearities from the head loss models for pipes, pumps, and valves respectively. The function $\m f^\mathrm{MV}(\m u)$ collects both the nonlinear head gain model of pumps $\m f^\mathrm{M}(\m u)$ in~\eqref{equ:head-flow-pump} and the nonlinear head loss model of valves $\m f^\mathrm{V} (\m u)$ in \eqref{equ:head-flow-valve}. 
Note that $\m h^{\mathrm{R}}$ is the constant heads at reservoirs. The DAE~\eqref{equ:dae-abstract} can be written as
\begingroup
\setlength\arraycolsep{0.3pt}
\begin{subequations}
	\begin{align*}
	\begin{bmatrix}
	\m O &  & & & \\
	& \m O & \ & & \\
	\  & \  & \m I & & \\
	\  & \  & \ & \m  O & \\
	\  & \  & \ & \ & \m O\\
	\end{bmatrix} \hspace{-4pt}\begin{bmatrix}
	\m x_1(k+1)\\
	\m x_2(k+1)\\
	\m x_3(k+1)\\
	\m v(k+1)\\
	\m u(k+1)
	\end{bmatrix}&= \begin{bmatrix}
	\m E_{x_1} & \m E_{x_2} & \m E_{x_3} & \m O& \m O\\
	\m F_{x_1} & \m F_{x_2}  & \m O & \m O& \m O\\
	\m O & \m O	&\m I & \m B_v & \m O \\
	\m O & \m O & \m O & \m G_{v} & \m G_{u} \\
	\m O & \m H_{x_2} & \m O& \m O & \m O 
	\end{bmatrix}\hspace{-4pt}
	\begin{bmatrix}
	\m x_1(k) \\
	\m x_2(k) \\
	\m x_3(k) \\
	\m v(k) \\
	\m u(k) \\
	\end{bmatrix}
	\\
	&+
	\begin{bmatrix}
	\m f^\mathrm{P}  \\
	\m f^\mathrm{MV} \\
	\m 0 \\
	\m 0  \\
	\m 0 
	\end{bmatrix}+\begin{bmatrix}
	\m O \,&\, \m O \\
	\m O \,&\, \m O \\
	\m O \,&\, \m O \\
	\m I \,&\, \m O \\
	\m O \,&\, \m I 
	\end{bmatrix} \begin{bmatrix}
	\m d(k) \\
	\m h^\mathrm{R}
	\end{bmatrix}.
	\end{align*}
\end{subequations}
\endgroup
We define  new vectors that collect all the variables as
$$ \m z(k) \triangleq \Bigl \lbrace \m x_1(k), \m x_2(k),\m x_3(k), \m v(k), \m u(k) \Bigr \rbrace,$$
$$ \m f( \m z(k)) \triangleq \Bigl \lbrace  \m f^\mathrm{P}, \m f^\mathrm{MV}  \Bigr \rbrace, \m l(k) \triangleq \Bigl \lbrace  	\m d(k), \m h^\mathrm{R} \Bigr \rbrace.$$
Hence, the compact state-space DAE model can be written as
\begingroup
\setlength\arraycolsep{0.3pt}
	\begin{equation}~\label{equ:descrwdn}
\boxed{	\m E_{z} \m z(k+1)   =   \m A_z 
	\m z(k) + \m B_{f} \m f(\m z(k)) + \m B_{l} \m l(k),}
	\end{equation}
	where state-space matrices $\m E_z, \m A_z, \m B_f$ and $\m B_l$ are obtained from~\eqref{equ:dae-abstract}. Next, we formulate the problem at hand.
\endgroup
\vspace{-0.1cm}
	\begin{rem}
	The discrete-time (DT) model~\eqref{equ:descrwdn} can also be written in continuous-time  (CT) by modifying some of the state-space matrices. The algebraic equations remain intact, seeing these algebraic constraints can be posed identically in CT or DT.  The only difference needed to obtain the state-space matrices in CT is via removing the $\m x_1(k)$ term (that is, the $h_i^{\mathrm{TK}}(k)$ term) and the discretization time-step $\Delta t$ from \eqref{equ:tankhead}.
	\end{rem}

\vspace{-0.4cm}

\section{Problem Formulation}~\label{sec:pbmform}

\vspace{-0.4cm}

The objective of this work is two-fold: to show that $\m f(\m z(k))$ is locally\footnote{Although we discuss local properties, the presented theory is valid \textit{semi-globally} in the sense that state variables are confined  within a compact state-space that encapsulates upper and lower bounds on, for example, flow rates. That is, the local assumption is global from the practical perspective. } Lipschitz and One-Sided Lipschitz (OSL)
 and to analytically compute the corresponding constants. The formal definitions of these two properties are as follows.
 \vspace{-0.1cm} 
\begin{mydef}~\label{def1}
	A function $\m f:\Omega\subseteq \mathbb{R}^n \rightarrow \mathbb{R}^m$ is Lipschitz continuous on $\Omega$ if there exists a constant $K(\Omega) \in \mathbb{R}_+$ such that
		\begin{align*}
			\normof{\m f (\m z_1) - \m f (\m z_2)}_2 \leq K \normof{\m z_1 - \m z_2}_2, \; \; \forall \; \m z_1, \m z_2 \in \Omega.
		\end{align*} 
A function $\m f$ is locally Lipschitz continuous if for every $\m y \in \Omega$ there exists an open neighborhood $U$ of $\m y$ such that the restriction $\m f |_U$ is Lipschitz continuous with constant $K(U)$.
\end{mydef}
\vspace{-0.3cm}
\begin{mydef}~\label{def2}
	A function $\m f:\Omega\subseteq \mathbb{R}^n \rightarrow \mathbb{R}^n$ is One-Sided Lipschitz on $\Omega$ if there exists a constant $L(\Omega) \in \mathbb{R}$ such that
		\begin{align*}
			\inp{\m f (\m z_1) - \m f (\m z_2)}{\m z_1-\m z_2} \leq L \normof{\m z_1 - \m z_2}_2^2, \; \; \forall \; \m z_1, \m z_2 \in \Omega.
		\end{align*} 
\end{mydef}
Definitions \ref{def1} and \ref{def2} quantitatively characterize the restricted growth of nonlinearities during WDN operation and can then be used for a variety of control and state estimation formulations~\cite{Yadegar2018,Phanomchoeng2010,ALESSANDRI2004,ichalal2012observer,zemouche2013lmi,Abbaszadeh2010,zhang2012full,benallouch2012observer}. Note also that a global property with respect to some domain $\Omega \subset \mathbb{R}^{n_m}$ is a \textit{semi-global} property with respect to $\mathbb{R}^{n_m}$.  The next section showcases analytical methods to compute the Lipschitz constant for $\m  f(\cdot)$.
\section{Computing Lipschitz Constants}\label{sec:Lipschitz} 
\indent We  show simple but effective methods to analytically compute the Lipschitz constant $K$ with respect to the flow rates. The domain is defined as the set of all possible flow rates reached during network operation. 
To prove $\m f(\m z(k))$ is locally Lipschitz, it suffices to show that $\m f^\mathrm{P}(\m v)$, $\m f^\mathrm{M}(\m u),$ and $\m f^\mathrm{V} (\m u)$ are Lipschitz with corresponding Lipschitz constants $K^\mathrm{P}, K^\mathrm{M},$ and $K^\mathrm{V}$. In particular, this section shows that
\begin{equation*}
{\normof{\m f (\m z_1) - \m f (\m z_2)}_2\leq K \normof{\m z_1 - \m z_2}_2 }, \; \; \forall \; \m z_1, \m z_2 \in \Omega,
\end{equation*}
in the domain of attainable controllable and uncontrollable flows. Constant $K \geq 0$ is the Lipschitz constant of the WDN under the assumption that all other parameters (such as pipe roughness coefficients and pipe parameters) are held constant. Since our approach is analytical, the scope of this work is restricted to the conditions under which parameters in the system remain differentiable. However, these assumptions coincide with those used in practice and are minimal, for a function that is Lipschitz continuous is necessarily differentiable almost everywhere. \\
\indent Throughout the paper we will utilize the following lemmas that are known in real analysis literature~\cite{federer2014geometric}. We reproduce their proofs and include them in Appendix~\ref{app:A}.
\vspace{-0.1cm}
\begin{mylem} \label{lem1}
	If $f:[a,b] \rightarrow \mathbb{R}$ is a Lipschitz continuous function with constant $K$ and differentiable then $K \coloneqq \sup_{[a,b]}|f'(x)|$.
\end{mylem}
\vspace{-0.3cm}
\begin{mylem}~\label{lem2}
	A function $\m f:\Omega \subseteq \mathbb{R}^n \rightarrow \mathbb{R}^m$ continuously differentiable on a convex set $\Omega$ is Lipschitz continuous if and only if $\normof{\m J_{\m f}(\m x) }_2$ is bounded where $\m J_{\m f}(\m x) $ is the Jacobian matrix of $\m f$. Furthermore,  $\sup_{\Omega} \normof{\m J_{\m f}(\m x) }_2 $ is a Lipschitz constant.
\end{mylem}
\vspace{-0.3cm}
\begin{rem}~\label{rem2}
From Lemma \ref{lem2} and the inequality $\normof{\cdot}_2 \leq \normof{\cdot}_F$ it follows immediately that
	\begin{align*}
		\normof{\m f (\m z_1) - \m f (\m z_2)}_2 \leq \sup_{\Omega} \normof{\m J_{\m f}(\m x) }_F \cdot \normof{\m z_1 - \m z_2}_2,
	\end{align*}
	where $\normof{\cdot}_F$ is the Frobenius norm. This result is used in \cite{nugroho2019tac} to estimate the upper bound of the optimal Lipschitz constant.
\end{rem}
\vspace{-0.1cm}
 Since our method requires finding individual $K^\mathrm{P,M,V}$ we first discuss finding $K^\mathrm{P}$. For simplicity, we  abuse the notation in each section by reusing variable notation instead of their full labels, such as $q$ for flow rates and $R$ for friction coefficients.  
 
 \vspace{-0.35cm}
 
 \subsection{Lipschitz  Constant for Headloss Models in Pipes}~\label{pipesec}
 Let $\m f^\mathrm{P}:\Omega^{\mathrm{P}} \subseteq \mathbb{R}^{n_p} \rightarrow \mathbb{R}^{n_p}$ be a function which maps the flow through pipes to a vector whose components characterize nonlinear head loss due to friction in pipes given in~\eqref{equ:pipeloss-abstracted}.
 \vspace{-0.1cm}
 \begin{asmp}~\label{asmp:pipe}
 	The $i^\mathrm{th}$ pipe has known a priori flow rate \\ $q_{i} \in [q_i^{\mathrm{min}}, q_i^{\mathrm{max}}]$ and fixed friction coefficient $R_i > 0$. Lower bound $q_i^{\mathrm{min}}$ can be negative, and parameter $\mu \geq 1$ in \eqref{equ:head-flow-pipe} is fixed and the same for all pipes.
 \end{asmp}
\vspace{-0.1cm}
\indent Assumption~\ref{asmp:pipe} is by no means restricting as typical values of $\mu$ are usually 1.852 or 2 according to Hazen-Williams, Darcy-Weisbach, or Chezy-Manning formulations~\cite{mays2010water}. In practice, the minimum and maximum pipe flow rates $[q_i^{\mathrm{min}}, q_i^{\mathrm{max}}]$ should not exceed the maximum flow rate provided by pumps. This can be obtained from the pump head loss curve~\eqref{equ:head-flow-pump} by setting $h^\mathrm{\mathrm{M}}_{ij}$ as zero. Hence, for each pump connecting nodes $i$ and $j$, the maximum flow rate is decided by
$$q_{ij}^{\mathrm{max}} = s_{ij} \left(\frac{h^\mathrm{s}}{r_{ij}}\right)^{\frac{1}{\nu}}.$$
\indent It is necessary that the domain of the reachable flow rates for the WDN, that is, the total Cartesian product of flow rate ranges for each pipe, is convex. However, since each pipe flow rate falls within a compact interval of Cartesian product, the domains of all pipes is automatically convex. In practice, this manifests as knowing the minimal and maximal flow rates of each pipe a priori which WDN operators know well from historical data. \\
\indent In the following proof we walk through the use of Lemma \ref{lem1} to show that each vector component of the nonlinear function $\m f^\mathrm{P}$ is Lipschitz, and then we use the results of Lemma \ref{lem2} to generalize for the entire vector function. In later proofs of this paper, we will skip to the generalization. The following proposition shows how $K^\mathrm{P}$ is computed.
\vspace{-0.1cm}
\begin{proposition}~\label{prop1}
	 Suppose that $\m f^\mathrm{P}:\Omega^{\mathrm{P}} \subseteq \mathbb{R}^{n_p} \rightarrow \mathbb{R}^{n_p}$ is the vector function of nonlinear terms describing the head loss in pipes for the entire WDN. With Assumption~\ref{asmp:pipe}, this function is Lipschitz with constant
	 \begin{equation}\label{equ:K^P}
	 \boxed{K^\mathrm{P} = \mu \max_i R_i(\max{\{|q_i^{\mathrm{min}}|, |q_i^{\mathrm{max}}|\}})^{\mu-1}.}
	 \end{equation}
\end{proposition}  
\vspace{-0.2cm}
\begin{proof}
 Suppose that $\m f^\mathrm{P}:\Omega^{\mathrm{P}} \subseteq \mathbb{R}^{n_p} \rightarrow \mathbb{R}^{n_p}$ where $\Omega^{\mathrm{P}}$ is the $n_p$ Cartesian product $[q_1^{\mathrm{min}}, q_1^{\mathrm{max}}] \times \dots \times [q_{n_p}^{\mathrm{min}}, q_{n_p}^{\mathrm{max}}]$. Each component $ \m f_i^{\mathrm{P}}(\m v)$ can be written as
	\begin{align*}
		 \m f^{\mathrm{P}}_i(\m v) = R_i{q_{i}}|q_{i}|^{\mu-1},
	\end{align*}
where $q_i \in [q_i^{\mathrm{min}}, q_i^{\mathrm{max}}]$ and $R_i, \mu$ are fixed. Let $\mathcal{I}_i$ denote the interval  $[q_i^{\mathrm{min}}, q_i^{\mathrm{max}}]$. For $\mu\geq 1$, $\m f^{\mathrm{P}}_i(\m v)$ is differentiable everywhere and since $\mathcal{I}_i$ is compact, by Lemma \ref{lem1} we obtain
	\begin{align*}
		 \left| \m f^{\mathrm{P}}_i(q_{i_1}) -\m f^{\mathrm{P}}_i(q_{i_2})\right| & \leq \sup_{\mathcal{I}_i}\left|\pder{\m f^{\mathrm{P}}_i}{q_{i}}\right|\cdot |q_{i_1} - q_{i_2}| \\
		 &= \sup_{\mathcal{I}_i} \mu R_i |q_{i}|^{\mu-1} \cdot |q_{i_1} - q_{i_2}| \\
		 &= K^\mathrm{P}_i |q_{i_1} - q_{i_2}|.
	\end{align*}
Next, since $\Omega^{\mathrm{P}}$ is convex, then we must show that the induced 2-norm of the Jacobian of $\m f^\mathrm{P}$ is bounded. For an arbitrary matrix $\m J$ the induced Euclidean norm is computed by	
	\begin{align*}
		\normof{\m J}_2 = \sigma^\mathrm{max}(\m J) =\sqrt{\lambda\m (\m J^{\top}\m J)},
	\end{align*}
where $\lambda(\m J^{\top}\m J)$ is the principal eigenvalue of $\m J^{\top}\m J$. Note that $\m f^{\mathrm{P}}_i (q_i)$ is not coupled with multiple components of $\m v$ such that the Jacobian matrix $\m J_{\m f^\mathrm{P}}(\m v)$ is diagonal. Since each $K^\mathrm{P}_i$ is finite, then each component of the Jacobian is necessarily bounded. Therefore, the norm is also bounded. By Lemma~\ref{lem2} we have
	\begin{align*}
		 \normof{\m f^\mathrm{P}(\m v_1) - \m f^\mathrm{P} (\m v_2)}_2 & \leq \sup_{\Omega^{\mathrm{P}}}\normof{\m J_{\m f^\mathrm{P}}(\m v)}_2 \cdot \normof{\m v_1 - \m v_2}_2.
	\end{align*}
Therefore, we finally get
	\begin{equation*}
		K^\mathrm{P} = \sup_{\Omega^{\mathrm{P}}}\normof{\m J_{f^\mathrm{P}}(\m v)}_2 = \sup_{\Omega^{\mathrm{P}}} \max_i \left|\pder{\m f^{\mathrm{P}}_i}{q_{i}}\right| = \max_i K^\mathrm{P}_i.
	\end{equation*}
	\vspace{-0.3cm}
\end{proof}


\subsection{Lipschitz Constant for Head-Gain Models through Pumps}
 To prove $\m f^{\mathrm{M}}( \m u)$ is locally Lipschitz we place restrictions on pump parameters typically used in literature. Note that a pump with a speed of zero simply reduces to a pipe so we assume that pump speeds are nonzero.
 \vspace{-0.1cm}
 	\begin{asmp} ~\label{asmp:pump}
	The $i^\mathrm{th}$ pump has flow rate $q_{i} \in [q_i^{\mathrm{min}}, q_i^{\mathrm{max}}]$ with $q_i^{\mathrm{min}} > 0$. The friction coefficient $r_i > 0$, pump speed \\ $s_i \in (0,1]$, and pump parameter $1\leq \nu_i \leq 3$ in~\eqref{equ:head-flow-pump} are fixed.
 	\end{asmp}
 \vspace{-0.1cm}
The head loss model for pumps is not necessarily defined for negative flow rates. Therefore, we assume that all pump flows do not change direction during operation. This is sensible as there would be minimal application during regular demand in operating a pump to slow down but not to stop and/or to reverse the natural flow rate. Note that Assumption~\ref{asmp:pump} is reasonable for actual networks we tested. For example, the value of $\nu_i$ is in $[1.1, 2.59]$ which can be corroborated through the open source tool EPANET and tens of WDN templates~\cite{rossman2000epanet}. Next, we show how $K^{\mathrm{M}}$ is computed.
\vspace{-0.1cm}
\begin{proposition}~\label{prop2}
	Let $\m f^{\mathrm{M}}:\Omega^{\mathrm{M}} \subseteq  \mathbb{R}^{n_m}  \rightarrow \mathbb{R}^{n_m}$ be the vector function of nonlinear terms describing the head loss in pumps.  With Assumption~\ref{asmp:pump}, the function $\m f^{\mathrm{M}}$ is Lipschitz with constant
		\begin{equation}\label{equ:K^M}
		\boxed{	K^\mathrm{M} =   \max_i \nu_i r_i (q_{i}^{\mathrm{max}})^{\nu_i-1} s_{i}^{2-\nu_i}.}
		\end{equation} 
\end{proposition}
\vspace{-0.2cm}
\begin{proof}
\indent Let $\m f^{\mathrm{M}}:\Omega^{\mathrm{M}} \subseteq \mathbb{R}^{n_m}  \rightarrow \mathbb{R}^{n_m}$ where $\Omega^{\mathrm{M}}$ is the $n_m$ Cartesian product of compact intervals $[q_i^{\mathrm{min}}, q_i^{\mathrm{max}}]$. The pump head gain for each component $i$ of $\m f^{\mathrm{M}}$ is given by
		\begin{align*}
			 \m f^{\mathrm{M}}_i(\m u) = -h_i^\mathrm{s}s_{i}^2 + r_i q_{i}^{\nu_i} s_{i}^{2-\nu_i},
		\end{align*}
where $q_i \in [q_i^\mathrm{min}, q_i^\mathrm{max}]$, and $s_i$, $\nu_i$, $r_i$, and $h_i^\mathrm{s}$ are each fixed in $\mathbb{R}_+$. Repeating the steps of the previous proof, by Lemma~\ref{lem2}  we have 
	\begin{align*}
		K^\mathrm{M} = \sup_{\Omega^{\mathrm{M}}}\normof{\m J_{\m f^{\mathrm{M}}}(\m u)}_2 = \sup_{\Omega^{\mathrm{M}}} \max_i \left|\pder{\m f^{\mathrm{M}}_i}{q_{i}}\right|.
	\end{align*}
This completes the proof.
\end{proof}
\vspace{-0.1cm}
Notice that if we assume that all pumps in the network have identical friction coefficient $r>0$, and letting
		$	\bar{q} = \max_i q_i^{\mathrm{max}},$
	then $K^\mathrm{M}$ can be obtained by mapping
		\begin{align*}
			g(s, q) \rightarrow
				 \begin{cases}
					g(s^{\mathrm{max}}, \bar{q}) & 1 \leq\nu \leq 2\\
					g(s^{\mathrm{min}}, \bar{q}) & \nu >2
				 \end{cases},
		\end{align*}
	where $s^{\mathrm{min}}/s^{\mathrm{max}}$ are the smallest/largest pump speeds in the WDN and
			$g(s,q) = \nu r q^{\nu-1} s^{2-\nu}.$
The next section investigates finding Lipschitz constants for the nonlinear valve head loss model. 
\vspace{-0.1cm}
\subsection{Lipschitz Constant for Valve Models}
The analysis for valves bears much similarity to that for pipes but with the addition of a valve openness parameter, which characterizes the occlusion of the flow. 
\vspace{-0.1cm}
	\begin{asmp}\label{asmp:valve} 
The $i^\mathrm{th}$ valve has known a priori flow rate $q_{i} \in [q_i^{\mathrm{min}}, q_i^{\mathrm{max}}]$, fixed friction coefficient $R_i >0$, and fixed valve openness $o_i \in (0,1]$. Parameter $1 \leq \mu \leq 3$ in \eqref{equ:head-flow-valve} is the same for all valves and fixed. 
	\end{asmp}
\vspace{-0.1cm}
In practice, this assumption holds true for a single time period. For extended periods, we choose the maximum of $o_i$. The computation of $K^\mathrm{V}$ is presented below.
\vspace{-0.1cm}
\begin{proposition}~\label{prop3}
	 Suppose that $\m f^\mathrm{V}:\Omega^{\mathrm{V}} \subseteq \mathbb{R}^{n_v} \rightarrow \mathbb{R}^{n_v}$ is the vector function of nonlinear terms describing the head loss in valves for the entire WDN. With Assumption~\ref{asmp:valve}, the function $\m f^\mathrm{V}$ is Lipschitz with constant
	 \begin{equation*}
	 \boxed{K^\mathrm{V} = \mu \max_i o_iR_i(\max{\{|q_i^{\mathrm{min}}|, |q_i^{\mathrm{max}}|\}})^{\mu-1}.}
	 \end{equation*}
\end{proposition}
\vspace{-0.2cm}
\begin{proof}
	The nonlinearities are described by $\m f^\mathrm{V} : \Omega^{\mathrm{V}} \subseteq \mathbb{R}^{n_v} \rightarrow \mathbb{R}^{n_v}$ and each component is given by
	\begin{align*}
		 \m f^{\mathrm{V}}_i (\m u) =  o_{i}R_i {q_{i}}|q_{i}|^{\mu-1}.
	\end{align*}
	This proof is similar to the proof of Proposition~\ref{prop1}, and the Lipschitz constant takes the form
	\begin{align*}
		K^\mathrm{V} = \sup_{\Omega^{\mathrm{V}}} \max_i \left|\pder{\m f^{\mathrm{V}}_i}{q_{i}}\right|.
	\end{align*}
	\vspace{-0.24cm}
\end{proof}
%
%
Given the three propositions, we now present an analytical way to compute Lipschitz constant $K$ for the entire network.
\vspace{-0.1cm}
\begin{theorem}~\label{thm1}
	If $K^\mathrm{P}, K^\mathrm{M}$ and $ K^\mathrm{V}$ are the Lipschitz constants for the pipes, pumps, and valves, respectively, then the Lipschitz constant $K$ of the WDN is given by
	\begin{equation}~\label{equ:lip}
		\boxed{K = \max{\{K^\mathrm{P}, K^\mathrm{M}, K^\mathrm{V}}\}.}
	\end{equation}
\end{theorem}
\vspace{-0.2cm}
\begin{proof}
	Since $\m f = \{\m f^\mathrm{P}, \m f^\mathrm{M}, \m f^\mathrm{V}\}$ and  the Jacobian of each of these components is diagonal, by Lemma \ref{lem2} we have that
\small		\begin{align*}
		 	K &= \sup_{\Omega} \nrm{\m J_{\m f}}_2 \\
			&{ = \max \left\{ \sup_{\Omega^{\mathrm{P}}} \max_i \left|\pder{\m f^{\mathrm{P}}_i}{q_{i}}\right| , \sup_{\Omega^{\mathrm{M}}} \max_j  \left|\pder{\m f^{\mathrm{M}}_j}{q_{j}}\right|, \sup_{\Omega^{\mathrm{V}}} \max_k \left|\pder{\m f^{\mathrm{V}}_k}{q_{k}}\right| \right\} }\\
			&= \max{\{K^\mathrm{P}, K^\mathrm{M}, K^\mathrm{V}}\}.
		\end{align*}
		\vspace{-0.2cm}
\end{proof}
In the subsequent section, we discuss the computation of OSL constant for WDNs.
\vspace{-0.3cm}

\section{Computing The One-Sided Lipschitz Constant}\label{sec:OSL} 
Various control-theoretic studies have shown Lipschitz based feedback control and state estimation algorithms can be conservative in terms of feasibility of the control/estimator designs via semidefinite programming (SDP) formulations~\cite{Abbaszadeh2010,benallouch2012observer,zhang2012full}, specifically for discrete-time models. This motivated the development of control/estimation methods for a less-restrictive assumption on the vector-valued nonlinearity, namely the one-sided Lipschitz (OSL).  Here, we show that
\begin{align*}
	\inp{\m f (\m z_1, \m c) - \m f (\m z_2, \m c)}{\m z_1-\m z_2}& \leq L \normof{\m z_1 - \m z_2}^2_2,
\end{align*}
where $\m c$ is the collection of fixed components of the network such as valve openness and pump speeds, and $L$ is the OSL we seek to compute. This computation can then be used as the SDP formulations in~\cite{Abbaszadeh2010,benallouch2012observer,zhang2012full} require this constant for, e.g., state estimator design.  Throughout the section, we consider that  all previous assumptions from the Lipschitz analysis section apply. 
Since the same assumptions carry over, the nonlinear head loss function $\m f$ is differentiable everywhere on the domain. Hence, the optimal OSL constant $L$ can be computed directly by taking the supremum of the logarithmic norm  $ \eta_2(\cdot)$ of the Jacobian matrix of $\m f(\cdot)$ as shown in \cite{abbaszadeh2010nonlinear}
\begin{align*}
	L = \sup_\Omega \eta_2(\m J_{\m f}) =  \sup_\Omega \lim_{\epsilon \rightarrow 0^+} \frac{\nrm{I + \epsilon \m J_{\m f} }_2-1}{\epsilon}.
\end{align*}
To that end, the following lemma computes the logarithmic norm for an arbitrary finite-dimensional real diagonal matrix, which is then used to compute the logarithmic norm $\eta_2(\m J_{\m f})$.
\vspace{-0.1cm}
\begin{mylem}~\label{lem3}
	If $\m D \in \mathbb{R}^{n \times n}$ is a diagonal matrix then
		\begin{align*}
			\eta_2(\m D) = \max_i \m D_{i},
		\end{align*}
		where
		\begin{align*}
			\eta_2(\m D) \coloneqq  \lim_{\epsilon \rightarrow 0^+} \frac{\nrm{\m I + \epsilon \m D}_2-1}{\epsilon}
		\end{align*}
	is the log norm with respect to the induced Euclidean norm.
\end{mylem}
\vspace{-0.2cm}
\begin{proof}
	Let $\m D \in \mathbb{R}^{n \times n}$ be a diagonal matrix. As mentioned in the proof of Proposition \ref{prop1} the induced 2-norm simplifies to
		\begin{align*}
			\eta_2(\m D) &=  \lim_{\epsilon \rightarrow 0^+} \frac{\sqrt{\lambda(\m I + 2\epsilon \m D + \epsilon^2 \m D^2)}-1}{\epsilon}
		\end{align*}
The principal eigenvalue $\lambda (\m I +2\epsilon \m D +\epsilon ^2 \m D^2)$ is the largest element in the diagonal entries of $\m I +2\epsilon \m D +\epsilon ^2 \m D^2$. When $\epsilon$ is small it is positive. Let $D_{i_0} = \max_i D_i$, then
		\begin{align*}
			\eta_2(\m D) &= \lim_{\epsilon \rightarrow 0^+} \frac{\sqrt{1 + 2\epsilon D_{i_0} + \epsilon^2 D_{i_0}^2}-1}{\epsilon}.
		\end{align*}
If $| D_1| > D_2$ and $D_1 < 0$, then
	\begin{align*}
	0 < \epsilon / 2 < \frac{D_2 -  D_1}{D_1^2-D_2^2}
	\end{align*}
implies $2\epsilon D_1 + \epsilon^2 D_1^2 < 2\epsilon D_2 + \epsilon^2 D_2^2$. It follows that
		\begin{align*} 
			\eta_2(\m D) &= \lim_{\epsilon \rightarrow 0^+} \frac{D_{i_0} + \epsilon D_{i_0}}{ \sqrt{(1 + 2\epsilon D_{i_0} + \epsilon^2 D_{i_0}^2)}} =  D_{i_0}.
		\end{align*}
		\vspace{-0.2cm}
\end{proof}
\vspace{-0.1cm}
We employ Lemma \ref{lem3} to simplify the computation of OSL constant $L$ for WDN in terms of Lipschitz constants for the pipes, pumps, and valves---summarized as follows.
\vspace{-0.1cm}
\begin{proposition}~\label{prop4}
 	The WDN One-Sided Lipschitz constant $L$ is always identical to the Lipschitz constant $K$,
 	\begin{equation*}
 		\boxed{L = \max{ \{K^\mathrm{P},K^\mathrm{M},K^\mathrm{V}\} }.}
 	\end{equation*}
 \end{proposition}
\vspace{-0.1cm}
 \begin{proof}
 	From the limit definition of the optimal OSL constant can be computed as
through $ 		L  =   \sup_\Omega \eta_2(\m J_{\m f_i} ).$
 	Since $\m J_{\m f}$ is diagonal, by Lemma~\ref{lem3} we have
 	\begin{align*}
 		L =  \sup_\Omega \max_i \m J_{\m f_i} .
 	\end{align*}
 From our analysis in Section \ref{sec:Lipschitz}, the non-zero components of $\m J_{\m f^\mathrm{P}}(\m v)$, $\m J_{\m f^\mathrm{M}}(\m u)$, and $\m J_{\m f^\mathrm{V}}(\m u)$ are always positive. Therefore,
 	 \begin{align*}
 		L &=  \sup_\Omega \max_i \m J_{\m f_i} =  \sup_\Omega \max_i |\m J_{\m f_i}| = K.
 	\end{align*}
\vspace{-0.2cm}
 \end{proof} 
 
 \section{Numerically Computing Lipschitz and OSL Constants}~\label{sec:global}
 Other than the analytical approach presented above, numerical computations of Lipschitz and OSL constants can also be performed via the utilization of \textit{(i)} interval-based algorithm or \textit{(ii)} point-based method. 
 
 For \textit{(i)}, the authors in \cite{nugroho2019tac} present an interval-based method to numerically estimate Lipschitz and OSL constants for arbitrary nonlinear vector-valued functions. The methods rely on the assumption that the set $\Omega$ is a compact, $n$-dimensional hyperrectangle and the function $\m f(\cdot)$ has bounded partial derivatives in $\Omega$. Note that this assumption is indeed satisfied in WDN, hence the interval-based approach is applicable to find Lipschitz and OSL constants for WDN. 
 The interval-based algorithm combines interval arithmetic with branch-and-bound (BnB) routines to find the smallest bounds that confine the largest value of a function. 
 Following \cite{nugroho2019tac}, the upper bound $\bar{K}$ for the actual Lipschitz constant $K$ can be computed as the solution to the following, constrained global maximization problem   
 \begin{align}
 \bar{K} = \sqrt{\max_{\m z\in \Omega}\sum_{i}\norm{\nabla \m f_i(\m{z})}_2^2}. \label{eq:gamma_Lipschitz}
 \end{align} 
 As mentioned in Remark \ref{rem2}, this is recognizable as the Frobenius norm of the Jacobian. The OSL constant $L$, as reproduced from literature such as \cite{abbaszadeh2010nonlinear}, can be computed as
 \begin{align}
 L =\lambda_{\mathrm{max}}\left(\dfrac{1}{2}\left(\m J_{\m f} + \m J^\top_{\m f}\right)\right), \label{eq:gamma_OSL}
 \end{align}     
 where $\m J_{\m f}$ denotes the Jacobian matrix of $\m f(\cdot)$. Since it is not obvious to compute the maximum eigenvalue of $\m J_{\m f}$, the authors in \cite{nugroho2019tac} provide theory that can be used to find upper bounds for $L$. Readers are referred to \cite{nugroho2019tac} for more details.
 
 	\begin{figure}
 	\centering
 	\includegraphics[width=0.75\linewidth]{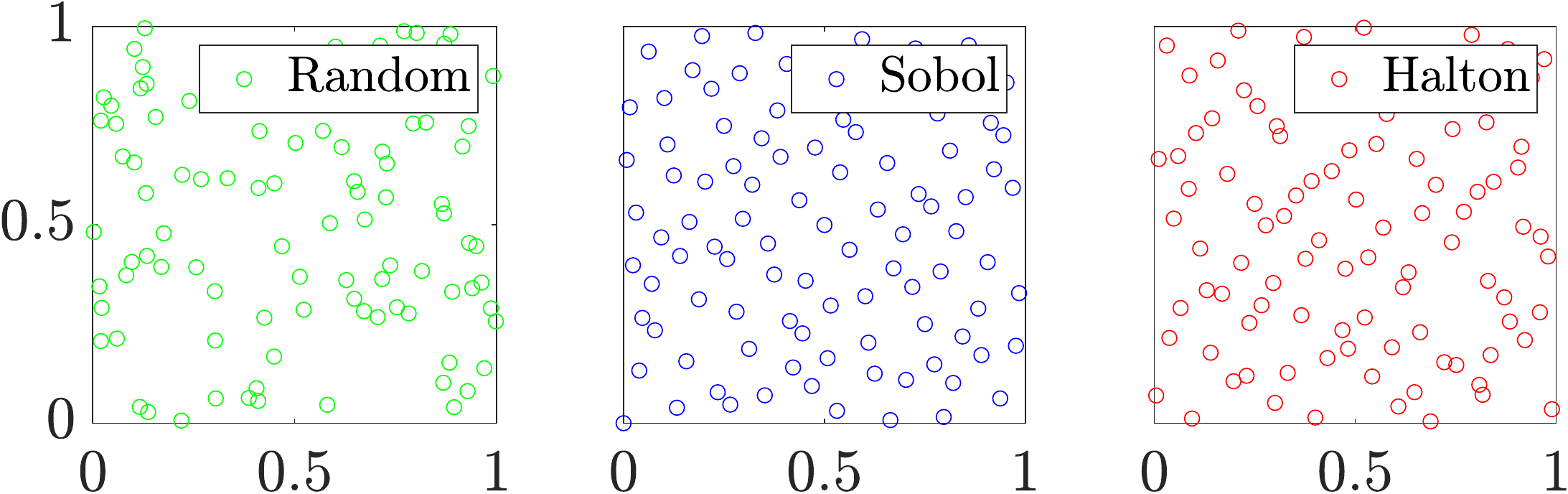}\vspace{-0.3cm}
 	\caption{Three different sampling methods. Note that the random points have more gaps and uneven distribution compared with Sobol and Halton.}
 	\vspace{-0.5cm}
 	\label{fig:sample3}
 \end{figure}
 
 The spatial BnB algorithm in the interval-based method can be succinctly explained as follows. In the branching step, the main problem---for instance, the maximization problem inside the square root term in \eqref{eq:gamma_Lipschitz}---is divided into smaller subproblems. This is performed by splitting $\Omega$. 
 Next, the corresponding upper and lower bounds of interval evaluation of each resulting subsets are computed accordingly. In the bounding step, the subsets not containing any maximizer are then removed. These two routines are performed iteratively until the algorithm terminates.
 In the context of Lipschitz constant, let $f^I_i(\cdot)$ be an interval extension of the objective function $\sum_i \norm {\nabla \m f_i(\cdot)}_2$ and define $C$ as a \textit{cover}---or a collection of subsets---in $\Omega$.
 In the BnB routines, it is crucial to ensure that all maximizers are always contained in $C$. The corresponding global upper and lower bounds are updated based on the interval evaluation of $f^I_i(\cdot)$ for each subset of $C$. The algorithm terminates when the \textit{optimality gap} (the distance between global upper and lower bounds) is sufficiently small.  Note that this method provides \textit{over-approximations} for Lipschitz and OSL constants.
 
Another approach to estimate Lipschitz and OSL constants is the \textit{(ii)} point-based method referenced above; see \cite{Nugroho2019Characterizing}. In principle, this method randomly samples the domain of interest $\Omega$ using a finite number of points and evaluates each point using the objective function to be maximized, ultimately finding the point which induces the largest value. As an illustration, in the case of Lipschitz constant, suppose that $ \mathcal{S}(\m z, s) = \left\{ \m z_j \right\}^s_{j = 1}$ represents the sequence of $s$ number of points with $\m z_j \in \Omega$ for each $j$. Then, we simply solve the following problem
\begin{align}
\ubar{K} =   \sqrt{\max_{\m z_j\in \mathcal{S}(\m z, s)}\sum_{i}\norm{\nabla \m f_i(\m{z}_j)}_2^2}, \label{eq:gamma_Lipschitz_under}
\end{align}
that is, finding $\m z_j$ from $\mathcal{S}(\m z, s)$ that maximizes \eqref{eq:gamma_Lipschitz_under}. The resulting maximum value is then regarded as the solution $\ubar{K}$. Realize that this approach involves evaluating the objective function of \eqref{eq:gamma_Lipschitz_under} as many as the number of points in the sequence. This method provides \textit{under-approximations} for Lipschitz and OSL constants.
 
 There exist several methods to construct $ \mathcal{S}(\m z, s)$. If the points are purely randomly generated, it is considered a \textit{Monte Carlo method} and the corresponding algorithm is termed \textit{pure random search}. 
 A more sophisticated method uses \textit{low-discrepancy sequences} (LDS) which are sequences of points with relatively small \textit{discrepancy} and are therefore distributed almost equally in the domain of interest.   
 The most important advantage of LDS---for example, Halton and Sobol sequences---is that it is guaranteed to converge to the optimal value as more points are used \cite{kucherenko2005application}. Fig. \ref{fig:sample3} depicts how purely random, Sobol, and Halton points are distributed inside a unit box. 
 It is worthwhile to mention that, however, as observed in the study \cite{Nugroho2019Characterizing}, the advantage of LDS over random sampling can be minuscule. That is, random sampling method can yield comparable results as those obtained from LDS---even slightly better in small cases.  

 	\begin{figure}[t]
 	\centering
 	\includegraphics[width=0.7\linewidth]{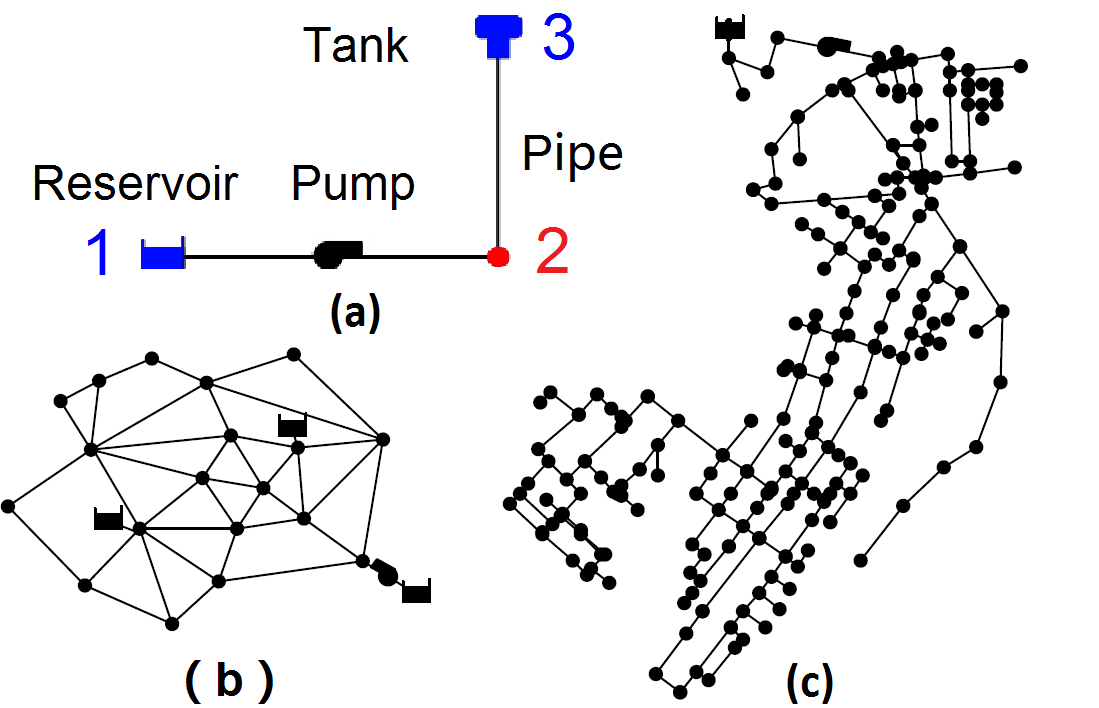}
 	\caption{(a) Three-Node network, (b) Anytown network, (c) OBCL network.}
 	\label{fig:scheme}
 	\vspace{-1.2em}
 \end{figure}

\section{Numerical Tests}~\label{sec:tests}
The aim of this section is to investigate how the analytical derivation of the Lipschitz constants (Section~\ref{sec:Lipschitz}) for various networks compare with the numerical computations (Section~\ref{sec:global}).  To that end, six WDN  templates (Three-Node~\cite{wang2020new}, Eight-node~\cite{rossman2000epanet}, Anytown~\cite{walski1987battle}, Net2, Net3, and OBCL networks~\cite{eliades2009epanet}) are used to assess our methods in finding Lipschitz constants $K$ via analytical and numerical methods. We do not compute OSL since $L$ has the same value as $K$; see Preposition~\ref{prop4}.
Fig.~\ref{fig:scheme} shows the network topology of the Three-Node, Anytown, and OBCL WDN---other test networks are not shown for brevity. The number of network's components is presented in the second column of Tab.~\ref{tab:my-table}. The Three-Node network shown in Fig.~\ref{fig:scheme}\textcolor{blue}{a} includes 1 junction, 1 pipe, 1 pump, and 1 reservoir. The pump parameters in~\eqref{equ:head-flow-pump} for this network are: $h^\mathrm{s}= 393.7008$, $r_{12} = 3.746 \times 10^{-6}$, and $\nu = 2.59$; the roughness parameter for Pipe 23 is $R_{23} = 
2.346 \times 10^{-6}$. 
The parameters for the other networks are utilized from~\cite{rossman2000epanet,walski1987battle,eliades2009epanet}. The bounds for the flow rates $q^\mathrm{min}$ and $q^\mathrm{max}$ for each link (pipe, pump, and valve) mentioned  in Section~\ref{sec:Lipschitz} are obtained from  EPANET after running the corresponding networks for a week with different water demand simulations. 
The numerical tests are performed with the help of the EPANET Matlab Toolkit~\cite{eliades2009epanet} on a Ubuntu 16.04 running an Intel(R) Xeon(R) CPU E5-1620 v3 @ 3.50GHz. All codes, parameters, tested networks, and results are available on Github~\cite{shenwangGP}.

	

	\begin{table}[t]
		\caption{Estimating optimal Lipschitz constant $K$ for all networks using three different methods.}
		\footnotesize
		\centering
		\label{tab:my-table}
		\setlength\tabcolsep{3pt}
		\renewcommand{\arraystretch}{1.5}
		\begin{tabular}{c|c|c|c|c|c}
			\hline
			\multirow{2}{*}{\textit{Networks}} & \multirow{2}{*}{\textit{\makecell{\# of com-\\ponents$^*$}} } & \multirow{2}{*}{\textit{Analytical$^\dagger$}} & \multicolumn{2}{c|}{\textit{\makecell{Point-based}}} & \multirow{2}{*}{\textit{\makecell{Interval- \\  based}}} \\ \cline{4-5}
			&  &  & \textit{max} & \textit{sqrt} &  \\ \hline
			\textit{\textit{\makecell{Three-Node\\network}}} & \makecell{\{1,1,1,\\1,1,0\}}  & 0.5023 & 0.5023 & 0.5023  & 0.5023 \\ \hline
			\textit{\textit{\makecell{Eight-node\\network}}} & \makecell{\{9,1,1,\\10,1,0\}} & 0.5850  & 0.5850 &  0.6084 & 0.6161 \\ \hline
			\textit{\textit{\makecell{Anytown\\network}}} & \makecell{\{19,3,0,\\40,1,0\}}   & 0.0903 & 0.0903 & 0.1250 & 0.1256  \\ \hline
			\textit{Net2} & \makecell{\{35,0,1,\\40,0,0\}}   & 0.0086 &  0.0086 &  0.0155& 0.0180 \\ \hline
			\textit{Net3} & \makecell{\{92,2,3,\\117,2,0\}}  & 0.0445 & 0.0445 &0.0689 & 0.0803 \\ \hline
			\textit{OBCL} & \makecell{\{262,1,0,\\288,1,0\}} & 0.9649  & 0.9649 & 1.0563 & 1.0817  \\ \hline \hline
			\multicolumn{6}{l}{\footnotesize{
					\makecell{$^*$\# of components: \{\# Junctions, \# Reservoirs, \# Tanks, \# Pipes, \# Pumps,\\ \# Valves\}.\\
				 Analytical$^\dagger$: this result is from Theorem \ref{thm1}, and specifically Equation \eqref{equ:lip}.  } }}
		\end{tabular}%
	\vspace{-0.4cm}
	\end{table}


As mentioned previously, this numerical test also considers point-based and interval-based methods to measure the conservativeness of analytical results. The interval-based method approximates Lipschitz constant $K$ via \eqref{eq:gamma_Lipschitz} whereas the point-based method uses two formulations: \eqref{equ:lip} (referred to as $\textit{max}$) and \eqref{eq:gamma_Lipschitz} (referred to as $\textit{sqrt}$). The \textit{max} mode can be used here since the induced Euclidean norm of the Jacobian coincides with the max norm for WDN.  Namely, the theoretical results of Lemma $\ref{lem2}$ are only valid for the induced Euclidean norm. 	 
The point-based method implements three different sampling methods (\textit{Random}, \textit{Sobol}, and \textit{Halton}) and uses various numbers of samples from $10^{1}$ to $10^{5}$ to show the relationship between accuracy (the quality of the approximation), sampling method, and numbers of samples.

The analytical and numerical results of this experiment are shown in Tab.~\ref{tab:my-table}. Taking the Three-Node network in Fig.~\ref{fig:scheme}\textcolor{blue}{a} as an example, the analytical Lipschitz constant for this network is $K = 0.5023$. This value is obtained from \eqref{equ:lip}.  Note that the Lipschitz constant for pump is $K^\mathrm{M} = 0.5023$, computed from~\eqref{equ:K^M}, while Lipschitz constant for pipe is $K^\mathrm{P} = 0.004$, computed from~\eqref{equ:K^P}, which is small due to the singular pipe's small roughness parameter $R_{23}$ and small flow rate. Due to Theorem \ref{thm1}, then from these results, the optimal analytical result is given as $K = \max\{K^\mathrm{P}, K^\mathrm{M}\} = 0.5023$. The estimated Lipschitz constant from interval-based method, computed via~\eqref{eq:gamma_Lipschitz}, gives an upper approximation of $\bar{K} = 0.5023$, which is very close with the analytical result.
Nonetheless, this agreement does not generalize to networks with pipes/valves that have larger parameters.
The point-based method, for the \textit{max} and \textit{sqrt} modes, surprisingly gives good under approximations in the simple Three-node network; the corresponding Lipschitz approximation is $\ubar{K} = 0.5023$, which is identical to the analytical one. However, this also does not generalize to other networks with larger number of pipes/valves with various parameters.

\begin{figure}[t]
		\centering
		\subfloat[\label{fig:sampling_b}]{\includegraphics[keepaspectratio=true,scale=0.7]{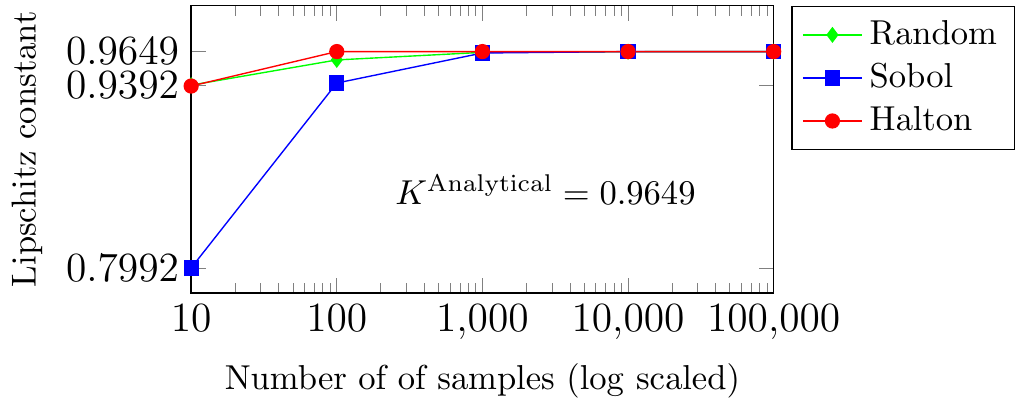}}{}\vspace{-0.2cm}\hspace{-0.0cm}\vspace{-0.2cm}
		\subfloat[\label{fig:sampling_c}]{\includegraphics[keepaspectratio=true,scale=0.7]{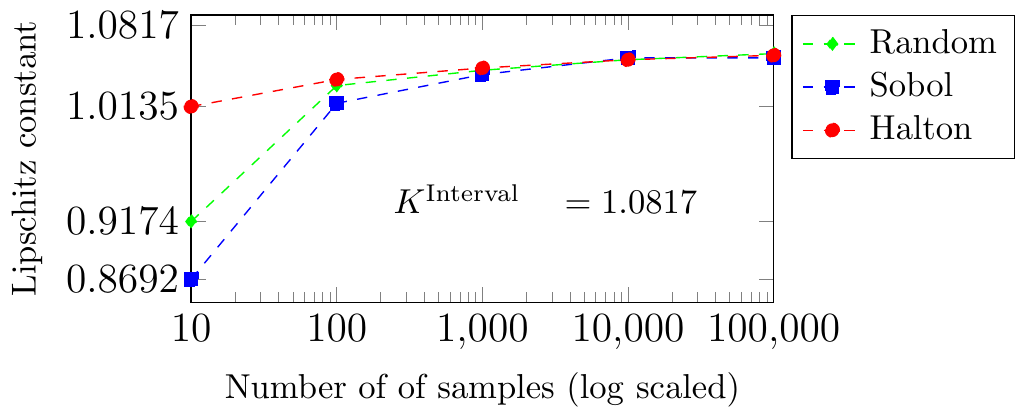}}{}\vspace{-0.2cm}\hspace{-0.0cm}\vspace{-0.05cm}
		\caption{Point-based method results for OBCL network under (a) \textit{max} mode  and (b) \textit{sqrt} mode.}
		\label{fig:sampling}\vspace{-0.30cm}
	\end{figure}

 
Next, we demonstrate the performance of point-based method in approximating Lipschitz constant for OBCL network shown in Fig.~\ref{fig:sampling}.
It can be seen from this figure that the point-based method under \textit{max} mode approaches the analytical Lipschitz constant $K^{\mathrm{Analytical}}$, while the solution obtained from \textit{sqrt} mode  approaches the one obtained from interval-based method $K^{\mathrm{Interval}}$.
We observe that, as the number of sampling point increases, the approximated Lispchitz constants from \textit{max} and \textit{sqrt} modes are approaching $K^{\mathrm{Analytical}}$ and $K^{\mathrm{Interval}}$  respectively, hence showing the increasing accuracy of point-based method. 
This phenomena is due to the property of LDS, where increasing number of samples produces better approximation. Interestingly, the results given by random sampling method stand between those from Halton and Sobol, indicating that LDS gives little advantage over pure random sampling.
	
\begin{figure}[t]
	\centering
	\includegraphics[width=0.7\linewidth]{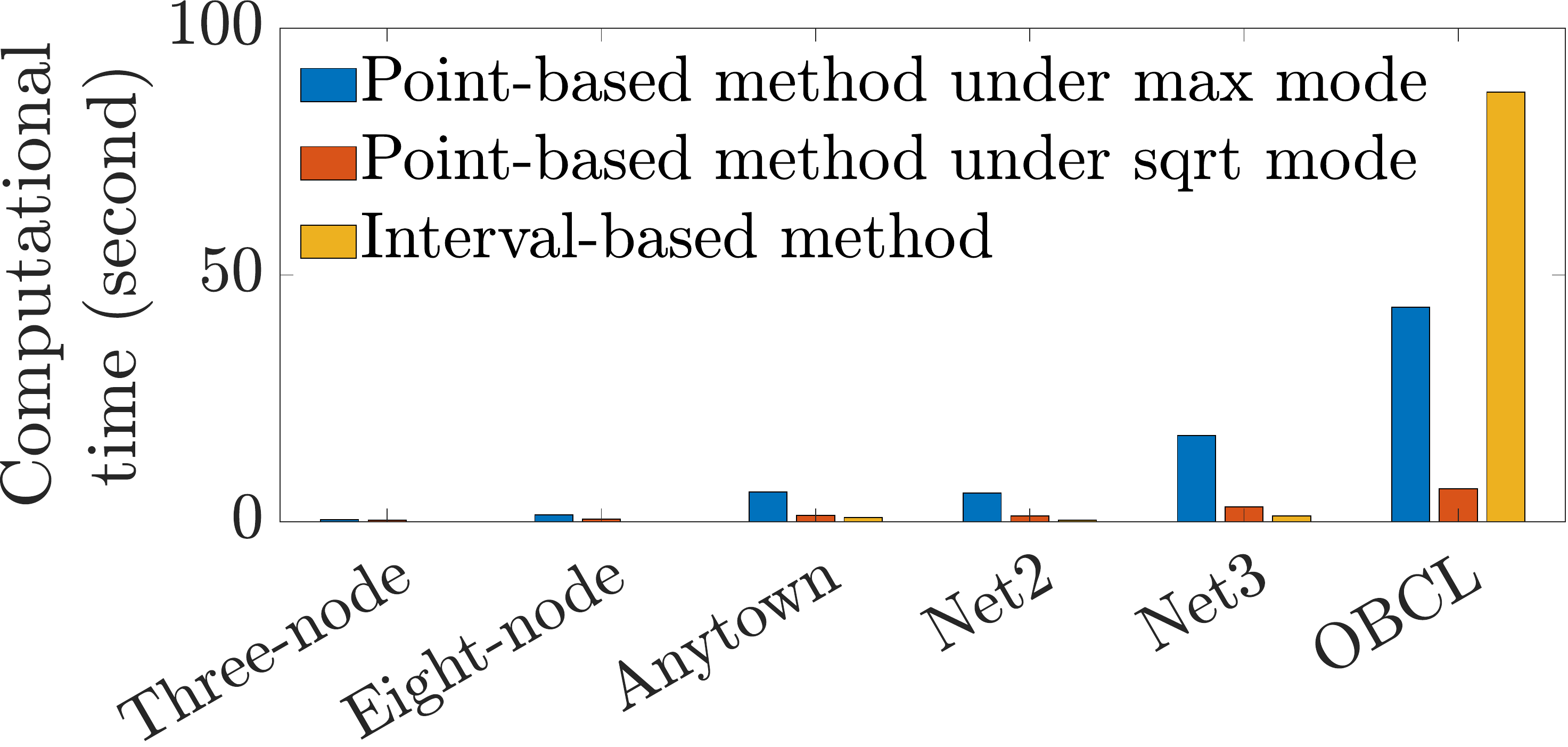}
	\caption{Computational time for the interval- and point-based methods.}
	\vspace{-0.4cm}
	\label{fig:computationaltime}
\end{figure}

Finally, we study the computational time of point-based and interval-based methods for each WDN network, which results are shown in Fig.~\ref{fig:computationaltime}.
In particular, for interval-based method, we consider the following optimality gap for the test networks: $1\times 10^{-2}$, $1\times 10^{-2}$, $1\times 10^{-5}$, $1\times 10^{-5}$, $2\times 10^{-3}$, and $8\times 10^{-2}$. The optimality gap defines the width of the interval bound for computing the Lipschitz constant: a larger gap requires less computational time for the same network, which justifies choosing a smaller one for the larger networks. 
In addition, to compensate the randomness in the point-based method, the displayed results are the average of running the simulation five times. 
The results suggest that the presented numerical methods are somewhat scalable even for a network with hundreds of nodes. 


Considering the derived analytical expressions and tested numerical algorithms for finding Lipschitz and one-sided Lipschitz constants, future work will focus on either applying or formulating state-feedback control and observer designs for the nonlinear DAE models of WDN hydraulics. Indeed, various studies~\cite{berger2018observers,koenig2006observer,wang2012observer,zulfiqar2016observer} have investigated this problem for a general Lipschitz or one-sided Lipschitz DAE model of dynamic systems.

\bibliographystyle{IEEEtran}
\bibliography{mybib}

\appendices 
\section{Proof of Lemmas~\ref{lem1} and~\ref{lem2}}\label{app:A}

\begin{proof}(of Lemma~\ref{lem1})
	Suppose that $f:[a,b] \rightarrow \mathbb{R}$ is Lipschitz continuous with constant $K$ and differentiable. Then for all $x,y \in [a,b]$, we have 
	\begin{align*}
\left |\frac{f(x)-f(y)}{x-y} \right | &\leq K.
	\end{align*}
	Taking the limit supremum of the inequality, we obtain
	\begin{align*}
	\limsup_{y\rightarrow x} \left |\frac{f(x)-f(y)}{x-y} \right |= K(x) \leq K.
	\end{align*}
	Since $[a,b]$ is compact the limit supremum is finite and coincides with the limit, and $f$ is differentiable on the interval $[a,b]$, we have
	\begin{align*}
	K(x) =\lim_{y\rightarrow x} \left| \frac{f(x)-f(y)}{x-y} \right | = |f'(x)|.
	\end{align*}
	Taking the supremum over all values in the compact interval we obtain a definition for $K$, that is,
	\begin{align*}
	K \coloneqq \sup_{[a,b]} K(x).
	\end{align*}
	\vspace{-0.7em}
\end{proof}

\begin{proof}(of Lemma~\ref{lem2})
	Suppose $\m f:\Omega \subseteq \mathbb{R}^n \rightarrow \mathbb{R}^m$ has Lipschitz constant $M$ and is continuously differentiable on $\Omega$. By the definition of the total derivative, for every vector $\m p$ with unit norm and $\epsilon > 0$, there exists $\delta > 0$ such that if $h < \delta$, then
	\begin{align*}
\left	| \frac{\normof{\m f(\m x + h\m p) - \m f(\m x)}_2}{h} - \normof{\m J_{\m f}(\m x)\m p}_2\, \right | < \epsilon.
	\end{align*}
	Hence,
	\begin{align*}
	\normof{\m J_{\m f}(\m x)\m p}_2 &< \frac{\normof{\m f(\m x + h\m p) - \m f(\m x)}_2}{h} + \epsilon < M + \epsilon.
	\end{align*}
	Since $\epsilon$ is arbitrary, we have for the induced 2-norm of $\m J_{\m f}(\m x)$
	\begin{align*}
	\normof{\m J_{\m f}(\m x)\m}_2 = \sup_{\normof{\m p} = 1} \normof{\m J_{\m f}(\m x)\m p}_2  \leq M.
	\end{align*}
	\indent Now suppose that $\normof{\m J_{\m f}(\m x) }_2$ is bounded by $M$ on a convex set $\Omega$. Then $\forall \m x,\m h \in \Omega$ and $t \in [0,1]$,	$\normof{\m J_{\m f}(\m x + t\m h)}_2 \leq M$. Therefore,
		\vspace{-0.7em}
	\begin{align*}
	\normof{\m f(\m x+\m h) - \m f(\m x)}_2 &= \left \normof{\int_{0}^{1} \m J_{\m f}(\m x + t\m h)\cdot \m h \, dt\right }_2  \\
	&\leq \int_{0}^{1}\left \normof{ \m J_{\m f}(\m x + t\m h)\cdot \m h\right }_2 \, dt \\
	&\leq M\normof{\m h}_2.
	\end{align*}
		\vspace{-2em}
\end{proof}

	\end{document}